\documentclass[12pt,a4paper]{article}

\usepackage{graphicx}
\usepackage{amsmath,amsfonts}
\usepackage{mathrsfs}
\usepackage[T1]{fontenc}
\usepackage[utf8]{inputenc}
\usepackage{authblk}
\usepackage{enumerate}
\usepackage{capt-of}
\usepackage{calc}
\usepackage{float}
\allowdisplaybreaks[1]

\def\Pb{\mathbf{P}}

\def\Ex{\mathbf{E}}

\def\AA{\mathbb{A}}
\def\BB{\mathbb{B}}

\def\UU{\mathbb{U}}

\def\sgn{{\rm sgn}} 
\def\1{\mbox{1\hspace{-.25em}I}}
\newcommand{\Liminf}{\mathop{\underline{\lim}}\limits}

\begin{document}
\title{Estimation of Cusp Location of Stochastic Processes: a Survey}
\author[1]{S.  Dachian}
\author[2]{N. Kordzakhia}
\author[3]{Yu.A. Kutoyants}
\author[4]{A. Novikov}
\affil[1]{\small University of Lille, Lille,  France}
\affil[2]{\small  Macquarie University, Sydney, Australia}
\affil[3]{\small  Le Mans University,  Le Mans,  France}
\affil[1,3]{\small National Research University ``MPEI'', Moscow, Russia}
\affil[3]{\small Tomsk State University, Tomsk, Russia}
\affil[4]{\small University of Technology Sydney, Australia}
\date{}

\maketitle
%{Running title : Estimation of Cusp }
\begin{abstract}
We present a review of some recent results on estimation of location parameter
for several models of observations with cusp-type singularity at the change
point. We suppose that the cusp-type models fit better to the real phenomena
described usually by change point models.  The list of models includes
Gaussian, inhomogeneous Poisson, ergodic diffusion processes, time series and
the classical case of i.i.d.\ observations. We describe the properties of the
maximum likelihood and Bayes estimators under some asymptotic assumptions. The
asymptotic efficiency of estimators are discussed as well and the results of
some numerical simulations are presented. We provide some heuristic arguments
which demonstrate the convergence of log-likelihood ratios in the models under
consideration to the fractional Brownian motion.

\end{abstract}
%\noindent MSC 2000 Classification: 62M02,  62G10, 62G20.

\bigskip
%\noindent 

\section{Models of observations}

We consider several models of observations having a cusp-type singularity at
the point of location parameter. Such models can be considered as a natural
alternative or extension to the change-point (in space and in time) models
with jumps in characteristics. The list of models of observations included in
this paper contains the case of independent identically distributed random
variables (i.i.d.  r.v.s, the change point in space), a signal in white
Gaussian noise (the change point in time), Poisson process (the change point
in time), the ergodic diffusion process (the change point in space) and
perturbed dynamical system (the change point in space).

 For general illustration of our approaches we start with the
model of signal  observed with a small additive white Gaussian noise:
\begin{equation}
\label{01}
{\rm d}X_t=S\left(\vartheta ,t\right){\rm d}t+\varepsilon {\rm d}W_t,\qquad
X_0=0,\quad 0\leq t\leq T .
\end{equation}
Here $S\left(\vartheta ,\cdot \right)\in {\cal L}_2\left[0,T\right]$ is the
signal, $\vartheta \in\Theta =\left(\alpha ,\beta \right)$ is an unknown 
parameter and $W_t$ is a standard Wiener process. The parameter $\varepsilon
>0$ is supposed to be small. We are interested in the following
  problem: how the asymptotic properties of  estimators depend on the order of regularity of the signal
  $S\left(\vartheta ,t\right)$? The special attention will be paid to the case of
  models with {\it cusp-type} singularity.

Note that the case of the observations $X^T=\left(X\left(t\right),0\leq
t\leq T=n\tau \right)$ of type  \eqref{01} with  $\tau $-periodic signal $S\left(\vartheta ,t\right)$
(the period $\tau $ is supposed to be known)
and $\varepsilon
=1$  can be reduced to the previous model  by setting
\[
X_t=\frac{1}{n}\sum_{j=1}^{n}\left[X\left(\tau
  \left(j-1\right)+t\right)-X\left(\tau \left(j-1\right)\right)\right],\qquad
0\leq t\leq \tau . 
\]
Now this averaging process satisfies \eqref{01} with $\varepsilon
=n^{-1/2}$ and another Wiener process (see also Section 1.3 for the case of
an inhomogeneous Poisson process with $\tau $-periodic intensity function).

To estimate the parameter 
$\vartheta $  and to
describe the asymptotic properties of the estimators as $\varepsilon
\rightarrow 0$ (a {\it small noise} asymptotics) we shall use the likelihood ratio function
\begin{equation}
\label{lr}
V\left(\vartheta ,X^T\right)=\exp\left\{\frac{1}{\varepsilon
  ^2}\int_{0}^{T}S\left(\vartheta ,t\right){\rm d}X_t-\frac{1}{2\varepsilon
  ^2}\int_{0}^{T}S\left(\vartheta ,t\right) ^2{\rm d}t\right\} ,\quad \vartheta
\in \Theta .
\end{equation}
In this paper we shall discuss  the asymptotic properties of the maximum likelihood
estimator (MLE) $\hat\vartheta _\varepsilon $ and Bayes   estimator (BE)
$\tilde\vartheta _\varepsilon 
$ for the quadratic loss
function. These estimators are defined by the relations
\[
V(\hat\vartheta_\varepsilon ,X^T)=\sup_{\vartheta \in \Theta
}V\left(\vartheta ,X^T\right),\qquad \tilde\vartheta _\varepsilon
=\frac{\int_{\alpha }^{\beta }\vartheta p\left(\vartheta
  \right)V\left(\vartheta ,X^T\right){\rm d}\vartheta }{\int_{\alpha }^{\beta
  } p\left(\vartheta \right)V\left(\vartheta ,X^T\right){\rm d}\vartheta}.
\]
We suppose that the    density $p\left(\cdot \right)$
of the prior distribution for $\vartheta $ is a 
positive continuous function. 

 It is known that if the signal $S\left(\vartheta ,t\right)$ is {\it smooth}
 w.r.t.\ $\vartheta $, then the MLE and BE are asymptotically normal and
 asymptotically efficient with the rate~$\varepsilon $, i.e. 
\[
\varepsilon ^{-1}(\hat\vartheta _\varepsilon -\vartheta )\Longrightarrow \zeta
,\qquad \varepsilon ^{-1}(\tilde\vartheta _\varepsilon -\vartheta
)\Longrightarrow \zeta ,\qquad \zeta \sim {\cal N}\left(0,{\rm
  I}\left(\vartheta \right)^{-1}\right),
\]
 where ${\rm I}\left(\vartheta \right)$ is the Fisher
 information \cite{IH1}. The arrow $\Longrightarrow $ means the convergence in
 distribution. 

In  contrast to the above smooth case, in the following classical change-point model of observations
\begin{equation}
\label{03}
{\rm d}X_t=\1_{\left\{t\geq \vartheta \right\}}{\rm
  d}t+\varepsilon {\rm d}W_t,\qquad X_0=0,\quad 0\leq t\leq T,
\end{equation}
the MLE $\hat\vartheta _\varepsilon$ and BE $\tilde\vartheta _\varepsilon  $ have  the rate of convergence $\varepsilon ^2$ (see  \cite{IH2}), i.e.
\[
\varepsilon ^{-2}\left(\hat\vartheta _\varepsilon
 -\vartheta \right)\Longrightarrow \hat\xi ,\qquad \varepsilon
 ^{-2}\left(\tilde\vartheta _\varepsilon 
 -\vartheta \right)\Longrightarrow \tilde\xi,
\]
where the random variables $ \hat\xi $ and $ \tilde\xi $  are defined by the relations
\begin{equation}
\label{04}
Z( \hat\xi)=\sup_{u\in {\cal R}}Z\left(u\right),\quad \qquad \tilde \xi
=\frac{\int_{-\infty }^{\infty }u Z\left(u\right){\rm d}u}{\int_{-\infty
  }^{\infty }Z\left(u\right){\rm 
    d}u}.
\end{equation}
Here $ Z\left(u\right)=\exp \{\gamma
W\left(u\right)-\frac{\left|u\right|}{2}\gamma^2\},\,u\in {\cal R} $,
$W\left(\cdot \right)$ is a two-sided Wiener process, $\gamma$ is a constant
(to be defined explicitly below). The uniqueness with probability $1$ of the
solution of the first equation was shown in~\cite{Pf82}.

Typically, real physical phenomena and technical devices have a transition
phase from one state to another which can be described in many ways, for
example, with use a smooth function (signal) having a very large Fisher
information. The important question is what happens if the regularity of the
signal is different of the supposed one, see some results in this direction in
\cite{Kut17}.  We consider here another model with the signals having {\it
  cusp}-type singularities. As alternative to \eqref{03} one can consider
\eqref{01} with, for example, the signal
\begin{equation}
\label{05}
S\left(\vartheta ,t\right)=\frac{1}{2}\left(1+{\sgn\left(t-\vartheta
  \right)}\left|\frac{t-\vartheta }{\delta }\right|^\kappa
\right)\1_{\left\{\left|t-\vartheta \right|\leq \delta \right\}}
+\1_{\left\{t>\vartheta + \delta \right\}} .
\end{equation}
Here $\delta >0$ is some small parameter and $\kappa \in (0,\frac{1}{2})$. We
suppose that \mbox{$\vartheta \in\left(\alpha ,\beta \right)$}, where $\alpha
>\delta $ and $\beta <T-\delta $.  This signal is a continuous function and
for small $\delta $ and $\kappa $ it can be a good ${\cal
  L}_2\left[0,T\right]$ approximation for the signal in~\eqref{01}. Note that
when the Fisher information does not exist, the problem of estimation
$\vartheta $ is singular.

\begin{figure}[H]
\centering
\includegraphics*[width=0.95\textwidth]{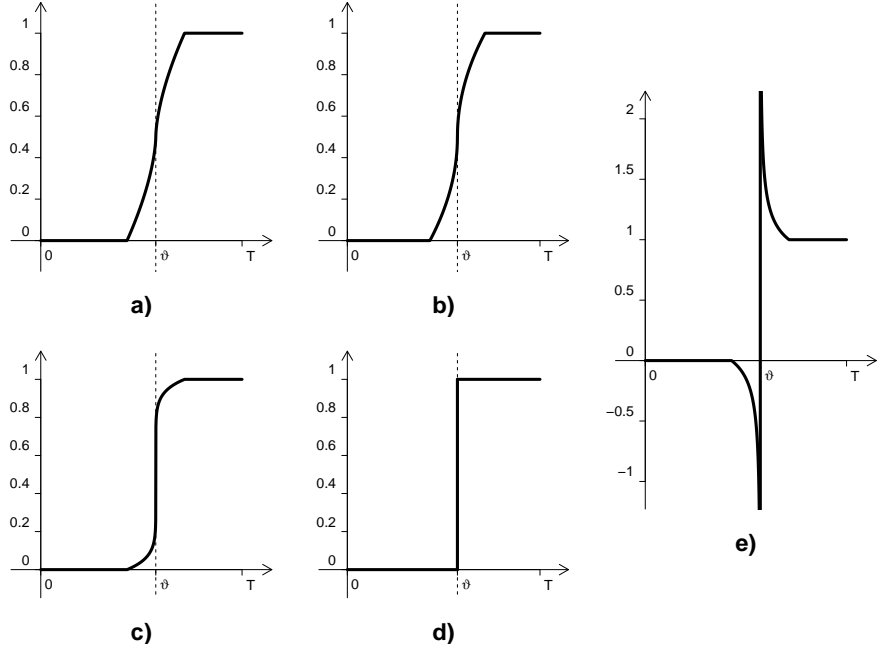}
\caption{\textbf{a)} $\kappa =\frac{5}{8}$, \textbf{b)} $\kappa =\frac{1}{2}$,
  \textbf{c)} $\kappa =\frac{1}{8}$, \textbf{d)} $\kappa =0$, \textbf{e)}
  $\kappa =-\frac{3}{8}$.}
\label{fig1}
\end{figure}

The examples of the corresponding curves of signals are given in
Figure~\ref{fig1}.

The asymptotic properties of the MLE and BE for $\vartheta $ under the
assumptions~\eqref{01} and \eqref{05} with $\kappa \in
(-\frac{1}{2},\frac{1}{2})$ are as follows:
\begin{equation}
\label{06}
\varepsilon ^{-\frac{2}{2\kappa +1}}\left(\hat \vartheta _\varepsilon
-\vartheta \right) \Longrightarrow \hat\xi  ,\qquad \varepsilon
^{-\frac{2}{2\kappa + 1}}\left(\tilde \vartheta _\varepsilon
-\vartheta \right) \Longrightarrow \tilde\xi ,
\end{equation}
where the random variables $\hat\xi$ and $\tilde\xi$ are defined in \eqref{04} with
\begin{equation}
\label{07}
Z\left(u\right)=\exp\left\{\gamma
W^H\left(u\right)-\frac{\left|u\right|^{2H}}{2}\gamma^2 \right\},\qquad
u\in{\cal R}.
\end{equation}
Here $H=\kappa +\frac{1}{2}$ is the {\it Hurst parameter}, 
\[
\gamma =\frac{1}{2\delta ^\kappa \Gamma_*} ,\qquad \Gamma_* ^2=\int_{-\infty
}^{\infty }\left[\sgn\left(s-1\right)\left|s-1\right|^\kappa
  -\sgn\left(s\right)\left|s\right|^\kappa\right]^2{\rm d}s
\]
 and $W^H\left(\cdot \right)$ is a  fractional
Brownian motion (fBm), i.e.\
 $W^H\left(u\right),u\in {\cal R}$,  is a 
Gaussian process,  $ W^H\left(0\right)=0, \Ex W^H\left(u\right)=0$ with the
covariance function
\begin{equation}
\label{fBm}
\Ex
W^H\left(u_1\right)W^H\left(u_2\right)=\frac{1}{2}\left(\left|u_1\right|^{2H}+
\left|u_1\right|^{2H}-\left|u_1-u_2\right|^{2H}\right).
\end{equation}
Note that it was first A.N.\ Kolmogorov \cite{Kol} who discussed this type of
Gaussian processes under the equivalent assumption to (\ref{fBm}):
\[
\Ex |W^H\left(u_1\right)-W^H\left(u_2\right)|^2=\left|u_1-u_2\right|^{2H}.
\]
The convergence in distribution, convergence of moments of these estimators
and the fact of the asymptotic efficiency of the BE were established in
\cite{CDK16} (for $\kappa \in \left(0,\frac{1}{2}\right)$) and in
\cite{KKNL17} (for $\kappa \in \left(-\frac{1}{2},0\right)$).

The rate of convergence of mean square error of the MLE and BE in the model
\eqref{01} with signal \eqref{05} essentially depends on $\kappa $.

We have the following asymptotics of the mean square error 

\begin{description}
\item[a)] {\it  Smooth signal, $\kappa >\frac{1}{2}$}: $\qquad \Ex_\vartheta \left(\hat \vartheta _\varepsilon -\vartheta \right) ^2\sim
\varepsilon ^2$ (see \cite{IH1}).

\item [b)] {\it Smooth signal with $\kappa =\frac{1}{2}$}: $\qquad \Ex_\vartheta \left(\hat \vartheta _\varepsilon -\vartheta \right) ^2\sim
\frac{\varepsilon^2 }{\ln \frac{1}{\varepsilon }}$.

\item [c)] {\it Continuous signal with  cusp $\kappa \in \left(0,\frac{1}{2}\right)$}:
  $\qquad \Ex_\vartheta \left(\hat \vartheta _\varepsilon -\vartheta \right)
  ^2\sim 
\varepsilon ^{\frac{4}{2\kappa +1}}$. {\it Note that 
  $2 <\frac{4}{2\kappa +1}< 4$} (see \cite{CDK16}).

\item [d)] {\it Discontinuous signal, $\kappa =0$}: $\qquad \Ex_\vartheta
  \left(\hat \vartheta _\varepsilon -\vartheta \right) ^2\sim \varepsilon ^4$
  (see \cite{IH2}).
\item [e)] {\it Discontinuous signal with  cusp $\kappa\in
  \left(-\frac{1}{2},0 \right)$}: $\qquad \Ex_\vartheta \left(\hat \vartheta _\varepsilon -\vartheta \right) ^2\sim
\varepsilon ^{\frac{4}{2\kappa +1}}$. {\it In this case $ \frac{4}{2\kappa +1}>4$}
  (see \cite{KKNL17}) .
\end{description}

It must be noted that if $\kappa \leq -\frac{1}{2}$, then the probability
measures generated by the observation $X^T$ are singular for all different
values $\vartheta$ and hence the parameter $\vartheta $ can be estimated
without error. To illustrate this fact, suppose that $S\left(\vartheta
,t\right)=\left|t-\vartheta \right|^\kappa $ and $\kappa \in
\left(-1,-\frac{1}{2}\right)$. Let $\kappa _*\in \left(-\frac{1}{2},0\right)$
such that $\kappa +\kappa _*\leq -1$ and the integral
\[
J\left(\vartheta \right)=\int_{0}^{T}\left|s-\vartheta \right|^{\kappa _*}{\rm
  d}X_s=\int_{0}^{T}\left|s-\vartheta \right|^{\kappa _*}\left|s-\vartheta_0
\right|^{\kappa }{\rm d}s+\varepsilon \int_{0}^{T}\left|s-\vartheta
\right|^{\kappa _*}{\rm d}W_s.
\]
It is easy to see that here the stochastic integral above is always finite with probability one and
the ordinary integral diverges at the point $\vartheta =\vartheta _0$. Using this
property we can construct an estimator of $\vartheta _0$ without error. For
example, $\vartheta _0$ is solution of the equation $J\left(\vartheta
\right)^{-1}=0 $. 

The case {\bf b)} can be treated similar  to {\bf a)}, using the same method
as in~\cite{IH1}.

The goal of this work is to present the review of the range of models with
cusp-type singularities and to describe the asymptotic properties of the MLE
and BE. Below we consider the i.i.d.\ random variables $X_1,\ldots,X_n$ with a
marginal density having the cusp-type singularity (as $n\rightarrow \infty $),
the ergodic diffusion processes $X_t,0\leq t\leq T$ with the trend coefficient
having cusp-type singularity (as $T\rightarrow \infty $), $\tau $-periodic
Poisson process \mbox{$X_t,0\leq t\leq n\tau $} with the intensity function
having a cusp-type singularity (as $n\rightarrow \infty $), the diffusion
processes with small diffusion coefficient $\varepsilon ^2$ and with the trend
coefficient having cusp-type singularity as $\varepsilon \rightarrow 0$. In
all such models the normalized likelihood ratio processes converge to the
process \eqref{07} with some constant~$\gamma $.  The proofs of week
convergence of the likelihood ratio processes to the exponent of fBm $W^{H}$
are rather tedious. Therefore, in this survey, we present only heuristic
arguments showing why convergence \eqref{07} should hold.  The detailed proofs
can be found in the cited papers.

 For illustrating our approaches to studying MLE and BE we start from
 \eqref{01} with the specified signal \eqref{05} under the assumption
 $\varepsilon \rightarrow 0$. The general technique to study for such models
 was developed by Ibragimov and Khasminskii in series of papers and can be
 found in their fundamental monograph \cite{IH81}. To use this technique we
 introduce the normalized likelihood ratio process
\[
Z_{\varepsilon }\left(u\right)=\frac{V\left(\vartheta _0 +\varphi _{\varepsilon }
  u,X\right)}{V\left(\vartheta _0 ,X\right)} ,\qquad u\in \UU_{\varepsilon
}=\left(\frac{\alpha -\vartheta _0 }{\varphi _{\varepsilon } },\frac{\beta
  -\vartheta _0 }{\varphi _{\varepsilon } } \right),
\]
where we denoted $\vartheta _0$ the true value and the normalizing function
$\varphi_\varepsilon =\varepsilon ^{\frac{1}{H}}$. In all problems under  consideration
in this paper we represent the log-likelihood ratio as follows
\[
\ln Z_\varepsilon \left(u\right)=\AA_\varepsilon
\left(u\right)-\BB_\varepsilon \left(u\right) +o\left(1\right)
\]
demonstrating the following weak convergences
\[
\AA_\varepsilon
\left(u\right) \Longrightarrow \gamma W^H\left(u\right) ,\qquad
\BB_\varepsilon \left(u\right)\Longrightarrow
\frac{\left|u\right|^{2H}}{2}\gamma ^2. 
\]
Here we would like to
show the universality of this local structure for rather different models of
observations without providing full technical details of the proofs.

 Suppose that we have already
proved the weak convergence  to the process $Z(u)$ in \eqref{07}
\[
Z_\varepsilon \left(\cdot \right)\Longrightarrow Z\left(\cdot \right).
\]
Then the limit distributions of the estimators can be  obtained as follows. For
the MLE $\hat \vartheta _\varepsilon $ we can write
\begin{align}
\label{08}
&\Pb_{\vartheta _0}\left(\frac{\hat\vartheta _\varepsilon-\vartheta_0
}{\varphi _\varepsilon }<x\right)=\Pb_{\vartheta _0}\left(\hat\vartheta
_\varepsilon< \vartheta_0+\varphi _\varepsilon x\right) \notag\\
&\qquad \qquad =\Pb_{\vartheta _0}\left\{ \sup_{\vartheta <\vartheta_0+\varphi
  _\varepsilon x} V\left(\vartheta ,X^T\right)>\sup_{\vartheta \geq
  \vartheta_0+\varphi _\varepsilon x} V\left(\vartheta ,X^T\right) \right\}
\notag\\
&\qquad \qquad =\Pb_{\vartheta _0}\left\{ \sup_{\vartheta <\vartheta_0+\varphi
  _\varepsilon x} \frac{V\left(\vartheta ,X^T\right)}{ V\left(\vartheta_0 ,X^T
  \right)}>\sup_{\vartheta \geq \vartheta_0+\varphi _\varepsilon x}
\frac{V\left(\vartheta ,X^T \right) }{V\left(\vartheta_0 ,X^T \right)}
\right\}\notag\\
&\qquad \qquad =\Pb_{\vartheta _0}\left\{ \sup_{u <x,u\in\UU_\varepsilon }
Z_\varepsilon \left(u\right)>\sup_{u \geq x,u\in\UU_\varepsilon }
Z_\varepsilon \left(u\right) \right\}\notag\\
&\qquad \qquad \longrightarrow \Pb_{\vartheta _0}\left\{ \sup_{u <x}
Z\left(u\right)>\sup_{u \geq x } Z\left(u\right) \right\} =\Pb_{\vartheta
  _0}\left(\hat \xi <x\right) .
\end{align}

For the BE $\tilde\vartheta _\varepsilon $ with the change of variable
$\vartheta_u=\vartheta _0+\varphi _\varepsilon u $ we have
\begin{align*}
\tilde\vartheta _\varepsilon &=\frac{\int_{}^{}\theta p\left(\theta
  \right)V\left(\theta ,X^T\right){\rm d}\theta }{\int_{}^{} p\left(\theta
  \right)V\left(\theta ,X^T\right){\rm d}\theta}=\vartheta _0+\varphi
_\varepsilon \frac{\int_{\UU_\varepsilon }^{}u p\left(\theta_u
  \right)V\left(\theta_u ,X^T\right){\rm d}u }{\int_{\UU_\varepsilon }^{}
  p\left(\theta_u \right)V\left(\theta_u ,X^T\right){\rm d}u}\\
&=\vartheta _0+\varphi _\varepsilon \frac{\int_{\UU_\varepsilon }^{}u
  p\left(\theta_u \right)\frac{V\left(\theta_u ,X^T\right)}{V\left(\vartheta
    _0,X^T\right)}{\rm d}u }{\int_{\UU_\varepsilon }^{} p\left(\theta_u
  \right)\frac{V\left(\theta_u ,X^T\right)}{V\left(\vartheta
    _0,X^T\right)}{\rm d}u}=\vartheta _0+\varphi _\varepsilon
\frac{\int_{\UU_\varepsilon }^{}u p\left(\theta_u \right)Z_\varepsilon
  \left(u\right){\rm d}u }{\int_{\UU_\varepsilon }^{} p\left(\theta_u
  \right)Z_\varepsilon \left(u\right){\rm d}u}.
\end{align*}
Hence
\begin{equation}
\label{09}
\frac{\tilde\vartheta _\varepsilon-\vartheta _0}{\varphi _\varepsilon
}=\frac{\int_{\UU_\varepsilon }^{}u p\left(\theta_u  
  \right)Z_\varepsilon \left(u\right){\rm d}u }{\int_{\UU_\varepsilon }^{}
  p\left(\theta_u 
  \right)Z_\varepsilon \left(u\right){\rm d}u}\Longrightarrow \frac{\int_{\cal
    R
  }^{}u Z\left(u\right){\rm d}u }{\int_{\cal R }^{} Z\left(u\right){\rm d}u}=\tilde
\xi.
\end{equation}
Thus we have (formally) verified \eqref{06}.

Moreover, in all such models under consideration the following lower bound for
the mean-square errors of all normalized estimators $\bar\vartheta
_\varepsilon $ holds:
\[
\lim_{\delta \rightarrow 0}\Liminf_{\varepsilon \rightarrow 0
}\sup_{\left|\vartheta -\vartheta _0\right|\leq \delta } \Ex_{\vartheta }
\left(\frac{\bar\vartheta _\varepsilon -\vartheta }{\varphi _\varepsilon
}\right)^2\geq \Ex_{\vartheta_0 } (\tilde \xi^2 )
\]
(see, e.g.\ \cite{CDK16}). This implies that an estimator $\vartheta
_\varepsilon ^*$ is asymptotically efficient (i.e. optimal) if for all
$\vartheta _0\in \Theta $ we have
\[
\lim_{\delta \rightarrow 0}\lim_{\varepsilon \rightarrow 0
}\sup_{\left|\vartheta -\vartheta _0\right|\leq \delta } \Ex_{\vartheta }
\left(\frac{\vartheta _\varepsilon ^*-\vartheta }{\varphi _\varepsilon
}\right)^2= \Ex_{\vartheta_0 } (\tilde \xi^2) .
\]

Note that if we put $\varphi _\varepsilon =\varepsilon
^{\frac{1}{H}} \gamma^{-\frac{1}{H}} $,  then
\[
Z_\varepsilon \left(\cdot \right)\Longrightarrow Z_0\left(\cdot \right),
\]
where $Z_0\left(\cdot \right)$ coincides with $Z\left(\cdot \right)$ in
\eqref{07} with $\gamma =1$.

Let $\hat \xi_0$ and $\tilde \xi_0$ be defined by \eqref{04} with
$Z\left(\cdot \right)$ replaced by $Z_0\left(\cdot \right)$. Obviously, the
distributions of $\hat \xi_0$ and $\tilde \xi_0$ do not depend on $\vartheta
_0$ and we obtain the following relations: $\hat \xi =\gamma ^{-\frac{1}{H}}
{\hat \xi_0}{}$ and $\tilde \xi = \gamma ^{-\frac{1}{H}}{\tilde \xi_0}{}$. In
particular, $\Ex_{\vartheta_0 } (\tilde \xi^2)=\gamma ^{-\frac{2}{H}}\Ex
(\tilde \xi_0^2) $.

It is of interest to compare the limit variances of $\hat\vartheta
_\varepsilon $ and $\tilde\vartheta _\varepsilon $ for the different values of
$H=\kappa+\frac{1}{2} $. In \cite{NKL14} it was shown via numerical
simulations that the limit values $\Ex (\hat\xi_0 ^2)$ could be essentially
larger than $\Ex (\tilde\xi_0 ^2)$. The results are presented in
Figure~\ref{fig2} for $H \in (0.4,1]$ or, correspondingly, $\kappa \in
  (-0.1,0.5] $. In Figure~\ref{fig3} we present the densities of the random
    variables $\hat\xi _0$ and $\tilde \xi _0$ obtained by the numerical
    simulations in \cite{KKNL17}. Note that on Panel B: $H=0.5$ the solid line
    shows the analytic curve for the density of MLE, this is the only case
    where the density is known in an analytic form, see \cite{KKNL17} for
    details.  \newdimen\myimageheight \newsavebox\myimageone
    \newsavebox\myimagetwo
    \sbox{\myimagetwo}{\includegraphics*[width=0.475\textwidth]{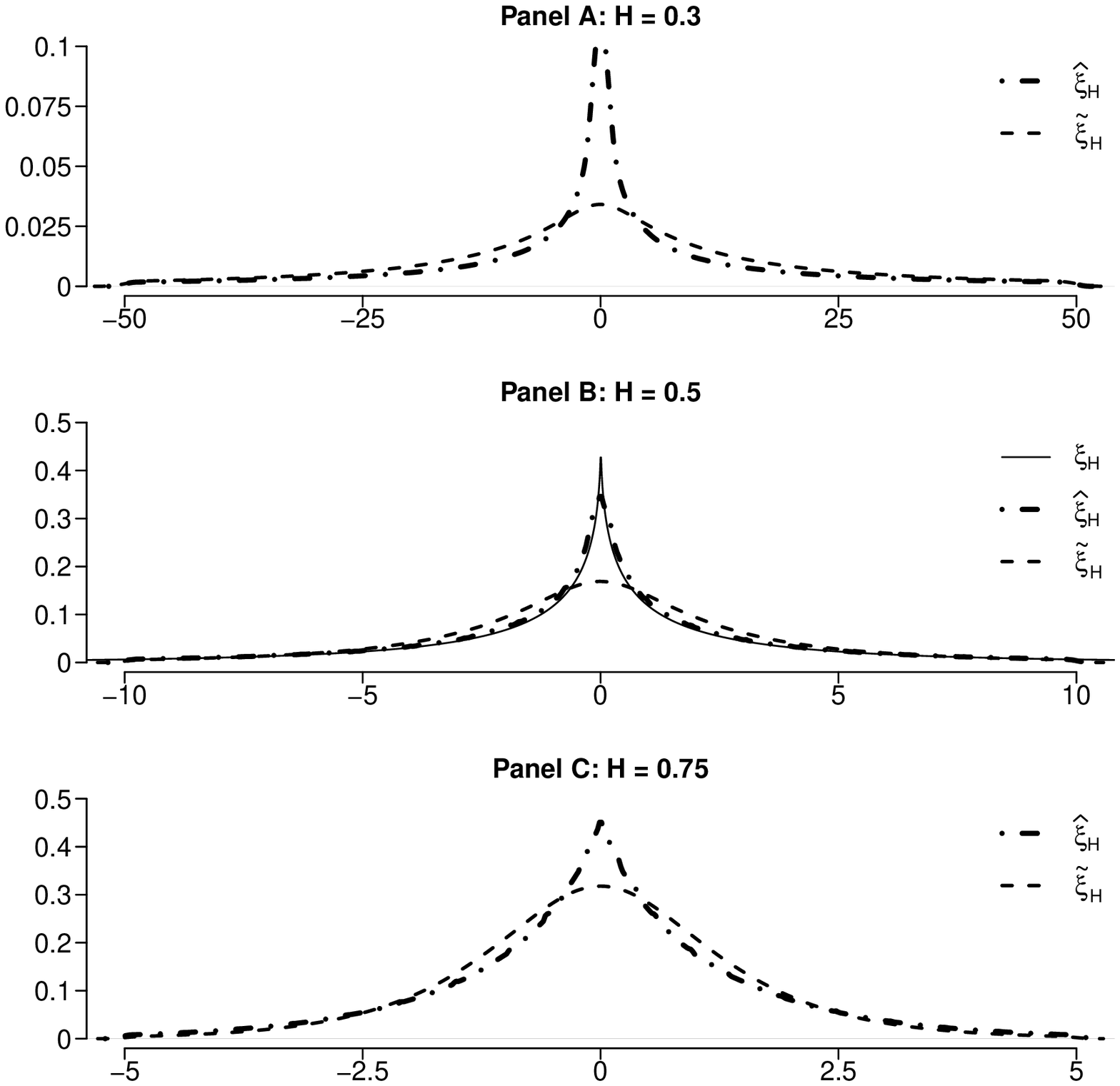}}
    \settototalheight{\myimageheight}{\usebox{\myimagetwo}}
    \sbox{\myimageone}{\includegraphics*[width=0.475\textwidth,
        totalheight=\myimageheight]{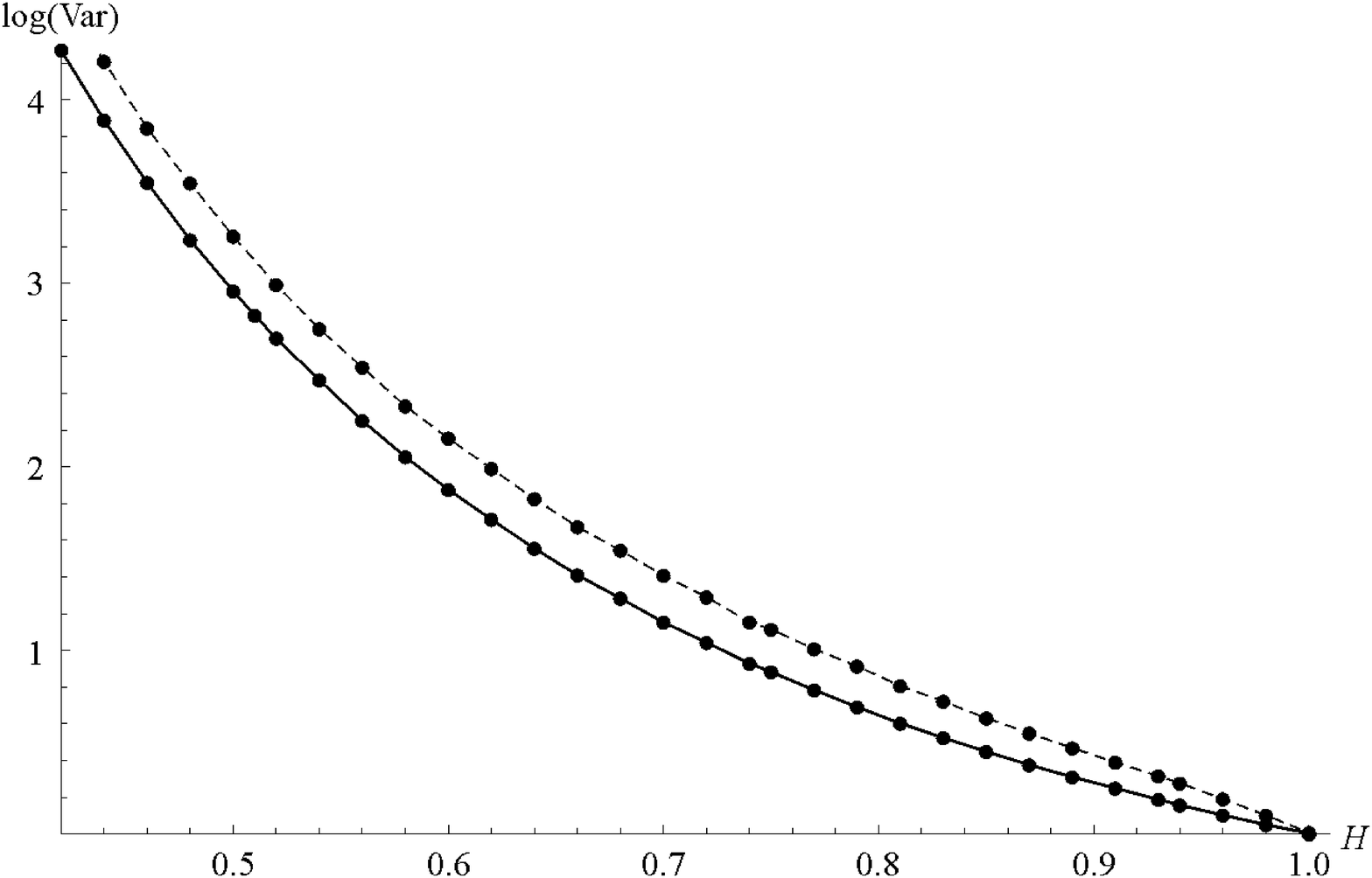}}

\begin{figure}[H]
\begin{minipage}{0.475\textwidth}
  \centering
  \usebox{\myimageone}
  \captionof{figure}{Limit curves of the values $\ln\Ex (\hat\xi_0 ^2)>\ln\Ex
    (\tilde\xi_0^2)$.}
  \label{fig2}
\end{minipage}
\hfill
\begin{minipage}{0.475\textwidth}
  \centering
  \usebox{\myimagetwo}
  \captionof{figure}{Densities of the MLE and BE for $H =0.3, 0.5, 0.75$.}
  \label{fig3}
\end{minipage}
\end{figure}

Below we consider several models of observations with the cusp-type
singularity which lead to the same limit likelihood ratio process
\eqref{07}. To demonstrate this we shall use the following representations of
the fBm
\[
W^H\left(u\right)=\Gamma_* ^{-1}\int_{-\infty }^{\infty
}\bigl[\sgn\left(v-u\right)\left|v-u\right|^\kappa
  -\sgn\left(v\right)\left|v\right|^\kappa \bigr] {\rm d}W\left(v\right).
\]
Here, obviously, $W^H\left(0\right)=0, \, \Ex W^H\left(u\right)=0$. Denoting
$s_{v,u}=\sgn\left(v-u\right)$ and $g_{v,u}=s_{v,u}\left|v-u\right|^\kappa -
s_{v,0}\left|v\right|^\kappa$ we can write
\begin{align*}
&\Ex W^H\left(u_1\right) W^H\left(u_2\right) =\Gamma_* ^{-2}\int_{-\infty
}^{\infty }g_{v,u_1}g_{v,u_2} {\rm d}v\\
&\quad =\frac{1}{2\Gamma_* ^{2}}\int_{-\infty }^{\infty }g_{v,u_1}^2 {\rm d}v
+\frac{1}{2\Gamma_* ^{2}}\int_{-\infty }^{\infty }g_{v,u_2}^2 {\rm d}v
-\frac{1}{2\Gamma_* ^{2}}\int_{-\infty }^{\infty }\bigl[g_{v,u_2}-g_{v,u_1}
  \bigr]^2 {\rm d}v\\
&\quad=\frac{1}{2}\left(\left|u_1\right|^{2H}+\left|u_2\right|^{2H}-\left|u_2-u_1\right|^{2H}
\right).
\end{align*}
Hence the process $W^H\left(\cdot \right)$ satisfies \eqref{fBm}. Here we used
the elementary identity $2ab=a^2+b^2-\left(a-b\right)^2$ and changed the
variables $v=s u_1$, $v=s u_2$, \mbox{$v=u_1+s\left(u_2-u_1\right)$}. 
 
We would like to mention here that a more
general representation of fBm is found in \cite{KKNL17}.

\subsection{Independent random variables} 

We suppose that  the observations 
$X^n=\left(X_1,\ldots,X_n\right)$ have the marginal density function
\[
f\left(x-\vartheta \right)=
{h\left(x-\vartheta\right)}\exp\left\{a\,\sgn\left(x-\vartheta
\right)\left|x-\vartheta \right|^\kappa \right\} ,\qquad \vartheta \in \Theta
=\left(\alpha ,\beta \right).
\]
The function $f\left(x \right)$ has at the point $x=0$ the cusp-type
singularity $\kappa \in \left(0,\frac{1}{2}\right)$.  The constant $a\not=0$
and the function $h\left(\cdot \right)$ are known, $h\left(0 \right)>0$.  We
suppose that the function $h\left(\cdot \right)$ is continuously
differentiable.  Our goal is to describe the behavior of the MLE
$\hat\vartheta _n$ and BE $\tilde\vartheta _n$. According to
\eqref{08}--\eqref{09} we describe the asymptotics of the normalized
likelihood ratio process
\[
Z_n\left(u\right)=\frac{V\left(\vartheta _0+\varphi
  _nu,X^n\right)}{V\left(\vartheta _0 ,X^n\right) },\qquad u\in
\UU_n=\left(\frac{\alpha -\vartheta _0 }{\varphi _n}, \frac{\beta  -\vartheta
  _0 }{\varphi _n}\right) ,
\]
where $\varphi _n=n^{-\frac{1}{2\kappa +1}}$. In particular, we have to verify
the convergence of $Z_n\left(u\right)$ to $Z\left(u\right)$ defined in
\eqref{07}. 
Let us see how the fBm appears in the limit log-likelihood ratio, using just
heuristic arguments. Below  we do it for $u\geq 0$. 
Let us denote $f_{u,x}= f\left(X_j-\vartheta _0-\varphi
  _nu\right)$, then we have
\begin{align*}
\ln Z_n\left(u\right)& =\sum_{j=1}^{n}\ln \frac{f\left(X_j-\vartheta
  _0-\varphi _nu\right)}{f\left(X_j-\vartheta _0\right)} =\sum_{j=1}^{n}\ln
\frac{f_{u,X_j}}{f_{0,X_j}}\\
&=\sum_{j=1}^{n}\ln \left(1+\frac{f_{u,X_j}-f_{0,X_j}}{f_{0,X_j}}\right)\\
&=\sum_{j=1}^{n}\frac{f_{u,X_j}-f_{0,X_j}}{f_{0,X_j}}-
\frac{1}{2}\sum_{j=1}^{n}\left(\frac{f_{u,X_j}-f_{0,X_j}}{f_{0,X_j}}\right)^2+o\left(1\right),
\end{align*}
where we used   Taylor expansion $ \ln \left(1+x\right)=x-\frac{x^2}{2}
+o\left(x^2\right)$. 

Introduce the further notations:
$g\left(u,x\right)=as_{x,u}\left|x-\vartheta_0- \varphi _nu\right|^\kappa$ and
$s_{x,u}=\sgn\left(x-\vartheta _0-\varphi _nu\right).  $ Using the expansion
$e^x-1=x+o\left(x\right)$ we can write
\[
\sum_{j=1}^{n}\left(\frac{f_{u,X_j}-f_{0,X_j}}{f_{0,X_j}}\right)^2=
\sum_{j=1}^{n}\left[g\left(u,X_j\right)-g\left(0,X_j\right)\right]^2\left(1+o\left(1\right)\right). 
\]
Recall the relation: for any continuous function $g\left(x\right)$ we have
\[
\frac{1}{n}\sum_{j=1}^{n}g\left(X_j\right)=\int_{-\infty }^{\infty
}g\left(x\right){\rm d}\hat F_n\left(x\right),\qquad \hat
F_n\left(x\right)=\frac{1}{n}\sum_{j=1}^{n} \1_{\left\{X_j<x\right\}},
\]
where $\hat F_n\left(x\right) $ is empirical distribution function. Let us
denote
\[
 \BB_n\left(u\right)=\frac12\sum_{j=1}^{n}{\left[g\left(u,X_j\right)-g\left(0,X_j\right)\right]^2}.
\]
Below $y=x-\vartheta _0, y=v\varphi _n, v=s u$:
\begin{align*}
\BB_n\left(u\right)&=n\,\frac{1}{2n}\sum_{j=1}^{n}\left[g\left(u,X_j\right)-g\left(0,X_j\right)\right]^2=\frac{n}{2}
\int_{ \cal R}^{}\left[g\left(u,x\right)-g\left(0,x\right)\right]^2{\rm d}\hat
F_n\left(x\right)\\
&\approx \frac{n}{2} \int_{ \cal
  R}^{}\left[g\left(u,x\right)-g\left(0,x\right)\right]^2{\rm d}
F\left(x-\vartheta _0\right)\\
&\approx \frac{a^2n}{2} \int_{ \cal R}^{}\left[\sgn\left(y-\varphi _nu\right)
  \left|y-\varphi _nu\right|^\kappa -\sgn\left(y\right) \left|y \right|^\kappa
  \right]^2f\left(y\right){\rm d} y \\
&=\frac{a^2n\varphi _n^{2\kappa +1}}{2} \int_{ \cal
  R}^{}\left[\sgn\left(v-u\right) \left|v-u \right|^\kappa -\sgn\left(v\right)
  \left|v \right|^\kappa \right]^2f\left(v\varphi _n\right){\rm d}v\\
&=\frac{a^2h\left(0\right)\left|u\right|^{2\kappa +1}}{2} \int_{ \cal
  R}^{}\left[\sgn\left(s-1\right) \left|s-1 \right|^\kappa -\sgn\left(s\right)
  \left|s \right|^\kappa \right]^2 {\rm d} s.
\end{align*}
Hence
\[
\BB_n\left(u\right)\longrightarrow \frac{a^2h\left(0\right)\left|u\right|^{2\kappa
    +1}}{2}\Gamma _*^2=\frac{\left|u\right|^{2\kappa +1}}{2} \gamma ^2.
\]
It is known that $B_n\left(x\right)=\sqrt{n}\left(\hat
F_n\left(x\right)-F_{\vartheta _0}\left(x\right)\right)\Rightarrow
B\left(F_{\vartheta _0}\left(x\right)\right) $, where $B\left(t\right)$, $t\in
\left[0,1\right]$ is a Brownian bridge, $B\left(t\right)=W\left(t\right)-t
W\left(1\right)$, $t\in \left[0,1\right]$.  Recall that $\Ex_{\vartheta
  _0}\frac{f_{u,X_j}-f_{0,X_j}}{f_{0,X_j}}=0 $.  Hence we can write formally
\begin{align*}
\AA_n\left(u\right)&=\sqrt{n}\int_{\cal
  R}^{}\frac{f_{u,x}-f_{0,x}}{f_{0,x}}{\rm d}B_n\left(x\right)\approx
\sqrt{n}\int_{\cal R}^{}\frac{f_{u,x}-f_{0,x}}{f_{0,x}}{\rm
  d}W_n\left(F_{\vartheta _0}\left(x\right)\right)\\
&\quad -\sqrt{n} W_n\left(1\right)\int_{\cal
  R}^{}\frac{f_{u,x}-f_{0,x}}{f_{0,x}}{\rm d}F_{\vartheta _0}\left(x\right)\\
&= \sqrt{n}\int_{\cal
  R}^{}\left[{g\left({u,x}\right)-g\left({0,x}\right)}\right]{\rm
  d}W_n\left(F_{\vartheta _0}\left(x\right)\right)+o\left(1\right).
\end{align*}

Below once more $y=x-\vartheta _0, y=v\varphi _n$ and $s_{y,u}=\sgn\left(y-u\right)$:
\begin{align*}
&\AA_n\left(u\right)=a\sqrt{n}\int_{\cal R}^{}\left[s_{y,u}\left|y-\varphi
  _nu\right|^\kappa -s_{y,0}\left|y\right|^\kappa\right]{\rm
  d}W_n\left(F\left(y\right)\right)+o\left(1\right) \\
&=a\sqrt{f\left(0\right)n}\varphi _n^{\kappa+\frac{1}{2} }\int_{\cal
  R}^{}\left[\sgn\left(v-u\right)\left|v-u\right|^\kappa
  -\sgn\left(v\right)\left|v\right|^\kappa\right] {\rm
  d}w_n\left(v\right)+o\left(1\right) \\
&\Longrightarrow a\sqrt{h\left(0\right)}\int_{\cal
  R}^{}\left[\sgn\left(v-u\right)\left|v-u\right|^\kappa
  -\sgn\left(v\right)\left|v\right|^\kappa\right]\; {\rm
  d}W\left(v\right)=\gamma W^H\left(u\right).
\end{align*}
Here we used the relations
\[
w_n\left(v\right)=\frac{W_n\left(F\left(\varphi _n
  v\right)\right)-W_n\left(F\left(0\right)\right)}{\sqrt{f\left(0\right)\varphi
    _n} }\Longrightarrow W\left(v\right) .
\]
 Therefore  
\[
Z_n\left(u\right)
=e^{\AA_n\left(u\right)-\BB_n\left(u\right)+o\left(1\right)}\Longrightarrow
Z\left(u\right)=e^{\gamma W^H\left(u\right)-\frac{\left|u\right|^{2\kappa
      +1}}{2} \gamma ^2}.
\]

The MLE $\hat\vartheta _n$ and the BE $\tilde \vartheta _n$ are consistent,
have different limit distributions 
\begin{align*}
&n^{\frac{1}{2\kappa +1}}\left(\hat\vartheta _n-\vartheta _0
\right)\Longrightarrow \hat\xi,\qquad Z(\hat\xi)=\sup_u Z\left(u\right) ,\\
& n^{\frac{1}{2\kappa +1}}\left(\tilde\vartheta _n-\vartheta _0
\right)\Longrightarrow \tilde\xi,\qquad \tilde\xi=\frac{\int_{}^{}u
  Z\left(u\right){\rm d}u}{\int_{}^{}Z\left(u\right){\rm d}u} ,
\end{align*}
their moments converge and the BE are asymptotically efficient:
\[
\lim_{\delta \rightarrow 0}\lim_{n \rightarrow \infty }\sup_{\left|\vartheta
  -\vartheta _0\right|\leq \delta }n^{\frac{2}{2\kappa +1}} \Ex_{\vartheta }
\left({\tilde\vartheta _n-\vartheta }\right)^2= \Ex_{\vartheta_0 } (\tilde
\xi^2).
\]

For the proof see \cite{IH81}, Section 6.4, where this case is called {\it
  singularity of order $2\kappa $ of the second type}. 

\textbf{Remark 1.} The so-called generalized  Gaussian distribution with the
density of the form
\[
f\left(x-\vartheta \right)=
{h\left(x-\vartheta\right)}\exp\left\{-a\left|x-\vartheta \right|^\kappa
\right\} ,\qquad \vartheta \in \Theta =\left(\alpha ,\beta \right), \kappa>0,
\]
was discussed by H.Daniels \cite{Daniels} for the case $\kappa >\frac12$. The
case $\kappa \in (0, \frac12)$ was first studied by Prakasa Rao \cite{PR68}.

\subsection{Signal in white Gaussian noise}

Consider the observations $X^T=\left(X_t,0\leq t\leq T\right)$ of the
stochastic process satisfying the equation
\begin{equation}
\label{2-0}
{\rm d}X_t=\left[a\,\sgn\left(t-\vartheta \right)\left|t-\vartheta
  \right|^k+h\left(t-\vartheta \right)\right]{\rm 
  d}t+\varepsilon {\rm d}W_t,\quad 
\quad 0\leq t\leq T .
\end{equation}
Here $X_0=0$, the parameters $a\not = 0, \kappa \in
\left(-\frac{1}{2},\frac{1}{2}\right)$ and function $h\left(\cdot \right)$ are
known and we have to estimate the parameter $\vartheta $. The function
$h\left(\cdot \right)$ is continuously differentiable and $\{W_t,0\leq t\leq
T\}$ is a standard Wiener process. The model defined by \eqref{01} and
\eqref{05} is similar to \eqref{2-0}.  We assume the $T$ is fixed and consider
the case $\varepsilon \rightarrow 0$ ({\it small noise asymptotics}).

  The unknown parameter $\vartheta \in \Theta =\left(\alpha ,\beta \right)$,
  where $0<\alpha<\beta <T$ and the likelihood ratio function is defined in
  \eqref{lr}, where the signal is of the form $S\left(\vartheta
  ,t\right)=a\,\sgn\left(t-\vartheta \right)\left|t-\vartheta
  \right|^k+h\left(t-\vartheta \right) $. We have to verify the convergence
\begin{equation}
\label{llr}
Z_\varepsilon \left(u\right)=\frac{V\left(\vartheta_0+\varphi _\varepsilon u
  ,X^T\right)  }{V\left(\vartheta_0 ,X^T\right)}\Longrightarrow 
Z\left(u\right)=\exp\left\{\gamma W^H\left(u\right)-\frac{\left|u\right|^{2H}}{2}\gamma ^2\right\}. 
\end{equation}
Here we set $\varphi _\varepsilon = {\varepsilon    ^{\frac{1}{\kappa
      +\frac{1}{2}}}}$ and $\gamma =a\Gamma_* $.

We shall demonstrate the validity of \eqref{llr}   for $u>0$. Below we denoted
$s_{t,u}=\sgn\left(t-\vartheta _0-\varphi _\varepsilon u\right)$. We have   
\begin{align*}
\ln Z_\varepsilon \left(u\right)&=\frac{a}{\varepsilon
}\int_{0}^{T}\left[s_{t,u}\left|t-\vartheta _0-\varphi _\varepsilon
  u\right|^\kappa -s_{t,0}\left|t-\vartheta _0\right|^\kappa\right] {\rm
  d}W_t\\
&\qquad -\frac{a^2}{2\varepsilon^2 }\int_{0}^{T}\left[s_{t,u}\left|t-\vartheta
  _0-\varphi _\varepsilon u\right|^\kappa -s_{t,0}\left|t-\vartheta
  _0\right|^\kappa\right]^2 {\rm d}t+o\left(1\right),
\end{align*}
because
\[
\frac{1}{\varepsilon^2 }\int_{0}^{T}\left[h\left(t-\vartheta _0-\varphi
  _\varepsilon u\right) -h\left(t-\vartheta _0\right)\right]^2 {\rm d}t\leq
C\frac{\varphi _\varepsilon ^2}{\varepsilon^2 } \leq C\,\varepsilon
^{\frac{2}{\kappa +\frac{1}{2}}-2}\longrightarrow 0.
\]
Further, 
\begin{align*}
&\BB_\varepsilon \left(u\right)=\frac{a^2}{2\varepsilon^2
  }\int_{0}^{T}\left[s_{t,u}\left|t-\vartheta _0-\varphi _\varepsilon
    u\right|^\kappa -s_{t,0}\left|t-\vartheta _0\right|^\kappa\right]^2 {\rm
    d}t\\
&\quad =\frac{a^2}{2\varepsilon^2 }\int_{-\vartheta _0}^{T-\vartheta
    _0}\left[\sgn\left(y-\varphi _\varepsilon u\right)\left|y-\varphi
    _\varepsilon u\right|^\kappa
    -\sgn\left(y\right)\left|y\right|^\kappa\right]^2 {\rm d}y\\
&\quad =\frac{a^2\varphi _\varepsilon ^{2\kappa +1}}{2\varepsilon^2
  }\int_{-\frac{\vartheta _0}{\varphi _\varepsilon }}^{\frac{T-\vartheta
      _0}{\varphi _\varepsilon
  }}\left[\sgn\left(v-u\right)\left|v-u\right|^\kappa
    -\sgn\left(v\right)\left|v\right|^\kappa\right]^2 {\rm d}v\\
&\quad =\frac{{a^2\left|u\right| ^{2\kappa +1}}}{2}\int_{-\frac{\vartheta
      _0}{u\varphi _\varepsilon }}^{\frac{T-\vartheta _0}{u\varphi
      _\varepsilon }}\left[\sgn\left(s-1\right)\left|s-1\right|^\kappa
    -\sgn\left(s\right)\left|s\right|^\kappa\right]^2 {\rm d}s\longrightarrow
  \frac{{\left|u\right| ^{2\kappa +1}}}{2}\gamma ^2,
\end{align*}
where we did the change of variables $t=y+\vartheta _0, y=\varphi _\varepsilon
v $  and $v=s u$. The similar change of variables in stochastic integral gives
us the following limit
\begin{align*}
\AA_\varepsilon \left(u\right)&=\frac{a}{\varepsilon
}\int_{0}^{T}\left[s_{t,u}\left|t-\vartheta _0-\varphi _\varepsilon
  u\right|^\kappa -s_{t,0}\left|t-\vartheta _0\right|^\kappa\right] {\rm
  d}W_t\\
&=\frac{a}{\varepsilon }\int_{-\vartheta _0}^{T-\vartheta
  _0}\left[\sgn\left(y-\varphi _\varepsilon u\right)\left|y-\varphi
  _\varepsilon u\right|^\kappa -\sgn\left(y\right)\left|y\right|^\kappa\right]
   {\rm d}\tilde W_y\\
&=\frac{a\varphi _\varepsilon ^{\kappa +\frac{1}{2}}}{\varepsilon
   }\int_{-\frac{\vartheta _0}{\varphi _\varepsilon }}^{\frac{T-\vartheta
       _0}{\varphi _\varepsilon }}\left[\sgn\left(v-u\right)\left|v-
     u\right|^\kappa -\sgn\left(v\right)\left|v\right|^\kappa\right] {\rm
     d}\tilde{\tilde W}_v\\
&={a}\int_{-\frac{\vartheta _0}{\varphi _\varepsilon }}^{\frac{T-\vartheta
       _0}{\varphi _\varepsilon }}\left[\sgn\left(v-u\right)\left|v-
     u\right|^\kappa -\sgn\left(v\right)\left|v\right|^\kappa\right] {\rm
     d}\tilde{\tilde W}_v\Longrightarrow a\Gamma_* W^H\left(u\right) .
\end{align*}
Therefore we verified the convergence \eqref{llr} and the MLE $\hat\vartheta
_\varepsilon $ and the BE~$\tilde\vartheta _\varepsilon $ are consistent, have
limit distributions 
The MLE $\hat\vartheta _n$ and the BE $\tilde \vartheta _n$ are consistent,
have different limit distributions 
\begin{align*}
&\varepsilon ^{-\frac{1}{\kappa +\frac{1}{2}}}\left(\hat\vartheta _\varepsilon
  -\vartheta _0 \right)\Longrightarrow \hat\xi,\qquad Z(\hat\xi)=\sup_u
  Z\left(u\right) ,\\
& \varepsilon ^{-\frac{1}{\kappa +\frac{1}{2}}}\left(\tilde\vartheta
  _\varepsilon-\vartheta _0 \right)\Longrightarrow \tilde\xi,\qquad
  \tilde\xi=\frac{\int_{}^{}u Z\left(u\right){\rm
      d}u}{\int_{}^{}Z\left(u\right){\rm d}u} ,
\end{align*}
their moments converge and the BE are asymptotically efficient.
For the proofs see \cite{CDK16}.

\subsection{Inhomogeneous Poisson processes} 

Suppose that we observe a trajectory of an inhomogeneous Poisson process
$X^T=\left(X_t,0\leq t\leq T\right)$ with $\tau $-periodic intensity function
$\lambda \left(t-\vartheta \right) $
admitting the  representation
\[
\lambda \left(t-\vartheta\right)=a\,\sgn\left(t-\vartheta \right)\left|t-\vartheta \right|^\kappa
+h\left(t-\vartheta \right),\qquad 0\leq t\leq \tau  ,
\]
on the first period and periodically continued on the whole real line.  Here
$\vartheta \in\Theta =\left(\alpha ,\beta \right)$, $0<\alpha <\beta <\tau
$. The function $\lambda \left(t\right)>0$, $t\in \left[0 , \tau \right]$ and
the parameter, $\kappa \in \left(0,\frac{1}{2}\right)$ are known. For
simplicity we assume that $T=n\tau $ and study asymptotics $n\rightarrow
\infty $.

The likelihood function is  (see \cite{LS01})
\[
V\left(\vartheta ,X^T\right)=\exp\left\{\int_{0}^{T}\ln \lambda
\left(t-\vartheta \right){\rm d}X_t-n\int_{0}^{\tau }\left[\lambda
  \left(t-\vartheta\right)-1\right]{\rm d}t \right\} ,\quad \vartheta \in\Theta .
\]
We have to show that the  normalized ($\varphi _n=n^{-\frac{1}{2\kappa +1}}$)
likelihood ratio process converges
\begin{equation}
\label{pl}
Z_n\left(u\right)=\frac{V\left(\vartheta_0+\varphi _nu ,X^T\right)
}{V\left(\vartheta_0 ,X^T\right)} \Longrightarrow
Z\left(u\right)=\exp\left\{\gamma
W^H\left(u\right)-\frac{\left|u\right|^{2H}}{2}\gamma ^2\right\}.  
\end{equation}
Here $\gamma =ah\left(0\right)^{-\frac{1}{2}}\Gamma_*$ with the same $\Gamma _*$ as before. 

Let us introduce the random processes
\begin{align*}
X_{j}\left(t\right)&=X_{\left\{\tau \left(j-1\right)+t\right\}}-X_{\left\{\tau
  \left(j-1\right)\right\}},\qquad 0\leq t\leq \tau, \quad j=1,\ldots,n,\\
W_n\left(t\right)&=\frac{1}{\sqrt{n}}\sum_{j=1}^{n}
\left[X_{j}\left(t\right)-\Lambda \left(\vartheta _0,t\right)\right] ,\quad
\Lambda \left(\vartheta _0,t\right)=\int_{0}^{t}\lambda \left(s-\vartheta
_0\right){\rm d}s ,
\end{align*}
and denote $\vartheta _u=\vartheta _0+\varphi _nu $, $s_{t,u}=\sgn
\left(t-\vartheta _u\right)$. Then 
using the relations $\ln \left(1+x\right)=x+ o\left(x\right),  x-\ln
\left(1+x\right)=\frac{x^2}{2} +o\left(x^2\right) $ and 
\begin{align*}
&\frac{\lambda \left(t-\vartheta _u\right)}{\lambda \left(t-\vartheta
    _0\right)}=1+\frac{\lambda \left(t-\vartheta _u\right)-\lambda
    \left(t-\vartheta _0\right)}{\lambda \left(t-\vartheta
    _0\right)}\\
&\qquad =1+\frac{as_{t,u}\left|t-\vartheta _u\right|^\kappa
    -as_{t,0}\left|t-\vartheta _0\right|^\kappa}{\lambda \left(t-\vartheta
    _0\right)}+\frac{h\left(t-\vartheta _u\right)-h\left(t-\vartheta
    _0\right)}{\lambda \left(t-\vartheta _0\right)}
\end{align*}
we can write
\begin{align*}
\ln Z_n\left(u\right)&= \sum_{j=1}^{n}\int_{0}^{\tau }\ln \left(\frac{\lambda
  \left(t-\vartheta _u\right)}{\lambda \left(t-\vartheta
  _0\right)}\right)\left[ {\rm d}X_{j}\left(t\right)-\lambda \left(t-\vartheta
  _0\right){\rm d}t\right] \\
&\qquad -n\int_{0}^{\tau }\left[\frac{\lambda \left(t-\vartheta
    _u\right)}{\lambda \left(t-\vartheta _0\right)}-1-\ln \left(\frac{\lambda
    \left(t-\vartheta _u\right)}{\lambda \left(t-\vartheta _0\right)}\right)
  \right]\lambda \left(t-\vartheta _0\right){\rm d}t\\
&= \sum_{j=1}^{n}\int_{0}^{\tau } \left(\frac{\lambda \left(t-\vartheta
  _u\right)-\lambda \left(t-\vartheta _0\right)}{\lambda \left(t-\vartheta
  _0\right)}\right)\left[ {\rm d}X_{j}\left(t\right)-\lambda \left(t-\vartheta
  _0\right){\rm d}t\right] \\
&\qquad -\frac{n}{2}\int_{0}^{\tau }\frac{\left(\lambda \left(t-\vartheta
  _u\right)-\lambda \left(t-\vartheta _0\right)\right)^2}{\lambda
  \left(t-\vartheta _0\right)}{\rm d}t+o\left(1\right)\\
&=a\sum_{j=1}^{n}\int_{0}^{\tau } \left(\frac{s_{t,u}\left|t-\vartheta
  _u\right|^\kappa -s_{t,0}\left|t-\vartheta
  _0\right|^\kappa}{a\,s_{t,0}\left|t-\vartheta
  _0\right|^\kappa+h\left(t-\vartheta _0\right)}\right)\left[ {\rm
    d}X_{j}\left(t\right)-\lambda \left(t-\vartheta _0\right){\rm d}t\right]\\
&\qquad -\frac{n a^2}{2}\int_{0}^{\tau }\frac{\left(s_{t,u}\left|t-\vartheta
  _u\right|^\kappa -s_{t,0}\left|t-\vartheta
  _0\right|^\kappa\right)^2}{a\,s_{t,0}\left|t-\vartheta
  _0\right|^\kappa+h\left(t-\vartheta _0\right)}\;{\rm d}t+o\left(1\right) ,
\end{align*}
because 
\[
n\int_{0}^{\tau }\frac{\left(h\left(t-\vartheta
  _u\right)-h\left(t-\vartheta _0\right)\right)^2}{\lambda \left(t-\vartheta
  _0\right)}{\rm d}t \leq C n\varphi _n^2=C\,n^{1-\frac{2}{2\kappa +1}}\longrightarrow 0.
\]
Below we change the variables $t=y+\vartheta _0, y=\varphi _n v, v=s u$
\begin{align*}
\BB_n\left(u\right)&=\frac{n a^2}{2}\int_{0}^{\tau
}\frac{\left(s_{t,u}\left|t-\vartheta _0-\varphi _nu\right|^\kappa
  -s_{t,0}\left|t-\vartheta
  _0\right|^\kappa\right)^2}{a\,s_{t,0}\left|t-\vartheta
  _0\right|^\kappa+h\left(t-\vartheta _0\right)}\;{\rm d}t\\
&=\frac{n a^2}{2}\int_{-\vartheta _0}^{\tau-\vartheta _0}
\frac{\left(\sgn\left(y-\varphi _nu\right)\left|y-\varphi _nu\right|^\kappa
  -\sgn\left(y\right)\left|y\right|^\kappa\right)^2}{a\,\sgn
  \left(y\right)\left|y\right|^\kappa+h\left(y\right)}\;{\rm d}y\\
&=\frac{n a^2\varphi _n^{2\kappa +1}}{2}\int_{-\frac{\vartheta _0}{\varphi
    _n}}^{\frac{\tau-\vartheta _0}{\varphi _n}}
\frac{\left(\sgn\left(v-u\right)\left|v-u\right|^\kappa
  -\sgn\left(v\right)\left|v\right|^\kappa\right)^2}{a\varphi _n^\kappa
  \sgn\left(v\right) \left|v\right|^\kappa+h\left(v\varphi _n^\kappa
  \right)}\;{\rm d}v\\
&=\frac{a^2\left|u\right|^{2\kappa
    +1}}{2h\left(0\right)}\int_{-\frac{\vartheta _0}{u\varphi
    _n}}^{\frac{\tau-\vartheta _0}{u\varphi _n}}
     {\left(\sgn\left(s-1\right)\left|s-1\right|^\kappa
       -\sgn\left(s\right)\left|s\right|^\kappa\right)^2}\;{\rm
       d}s+o\left(1\right)\\
&\quad \longrightarrow \frac{a^2\left|u\right|^{2\kappa
         +1}}{2h\left(0\right)}\Gamma_* ^2=\frac{\left|u\right|^{2\kappa
         +1}}{2}\;\gamma ^2 .
\end{align*}
The same change of variables ($t=\vartheta _0+\varphi _n v$) in the stochastic
integral provide us the relations
\begin{align*}
&\AA_n\left(u\right)=a\sqrt{n}\int_{0}^{\tau }
  \left(\frac{s_{t,u}\left|t-\vartheta _0-\varphi _nu\right|^\kappa
    -s_{t,0}\left|t-\vartheta _0\right|^\kappa}{a\,s_{t,0}\left|t-\vartheta
    _0\right|^\kappa+h\left(t-\vartheta _0\right)}\right) {\rm
    d}W_{n}\left(t\right)\\
&\quad =\frac{a\sqrt{n}\varphi _n^{\kappa
      +\frac{1}{2}}}{\sqrt{h\left(0\right)}}\int_{-\frac{\vartheta _0}{\varphi
      _n}}^{\frac{\tau -\vartheta _0}{\varphi _n} }
  \left({\sgn\left(v-u\right)\left|v-u\right|^\kappa
    -\sgn\left(v\right)\left|v\right|^\kappa}\right) {\rm
    d}w_n\left(v\right)\\
&\quad \Longrightarrow \frac{a}{\sqrt{h\left(0\right)}}\int_{-\infty }^{\infty
  }\left({\sgn\left(v-u\right)\left|v-u\right|^\kappa
    -\sgn\left(v\right)\left|v\right|^\kappa}\right) {\rm d}W\left(v\right)=
  \gamma \;W^H\left(u\right),
\end{align*}
where
\[
w_n\left(v\right)=\frac{W_{n}\left(\vartheta _0+\varphi _n
  v\right)-W_{n}\left(\vartheta _0 \right)}{\sqrt{\varphi
    _n}}\Longrightarrow  W\left(v\right). 
\]
Therefore we have the convergence \eqref{pl}.  The MLE $\hat\vartheta _n$ and the BE
$\tilde\vartheta _n$ are consistent, have limit distributions 
\[
n^{\frac{1}{2\kappa +1}}{{\left(\hat\vartheta _n -\vartheta
    _0\right)}}\Longrightarrow \hat \xi ,\qquad \qquad \qquad
n^{\frac{1}{2\kappa +1}}{{\left(\tilde\vartheta _n -\vartheta
    _0\right)}}\Longrightarrow \tilde \xi ,
\]
the moments of these estimators converge and the BE are asymptotically
efficient.  
For the proofs see \cite{D03}.

\subsection{Ergodic diffusion process} 

Consider the observations $X^T=\left(X_t,0\leq t\leq T\right)$ of the ergodic
diffusion process
\begin{equation}
\label{4-0}
{\rm d}X_t=\left[a\,\sgn\left(X_t-\vartheta \right)\left|X_t-\vartheta
  \right|^\kappa +h\left(X_t-\vartheta \right)\right]{\rm d}t + {\rm d}W_t,\quad 
X_0.
\end{equation}
Here $
a\not=0$, $\kappa \in \left(0,\frac{1}{2}\right)$, the function $h\left(\cdot
\right)$ is  known and has bounded derivative. 
Moreover we suppose that the function 
$S\left(x\right)=a\,\sgn\left(x\right)\left|x
  \right|^\kappa +h\left(x \right)$ is such that the conditions ${\cal ES}$, ${\cal EM}$ and
${\cal A}_0\left(\Theta \right)$ in \cite{Kut04} are fulfilled. For example,
these conditions are fulfilled if $h\left(x\right)=-bx$ with $b>0$. These
conditions provide the 
existence and uniqueness  of the solution of this equation, the existence of
finite invariant measure with the density
\[
f\left(\vartheta,x \right)=f\left(x-\vartheta \right)=G\exp\left\{2\int_{0}^{x-\vartheta
}\left[a\,\sgn\left(z\right)\left|z\right|^\kappa +h\left(z\right)\right]{\rm d}z\right\}
\]
and finiteness of all polynomial moments. Here $G>0$ is the normalizing constant.
The likelihood function is (see \cite{LS01})
\[
V\left(\vartheta ,X^T\right)=\exp\left\{\int_{0}^{T}{S\left(X_t-\vartheta\right)}{\rm 
  d}X_t-\frac{1}{2}\int_{0}^{T}{S\left(X_t-\vartheta\right)^2}{\rm
  d}t\right\},\qquad \vartheta \in \Theta . 
\]
 For the normalized likelihood ratio process we have to show the convergence
\begin{equation}
\label{4-1}
Z_T\left(u\right)=\frac{V\left(\vartheta_0+\varphi _T u ,X^T\right)
}{V\left(\vartheta_0 ,X^T\right)}\Longrightarrow
Z\left(u\right)=\exp\left\{\gamma
W^H\left(u\right)-\frac{\left|u\right|^{2H}}{2}\gamma ^2\right\}.
\end{equation}
Here $\varphi _T=T^{-\frac{1}{2\kappa +1}}$ and $\gamma =a\Gamma_* G^{1/2} $. 

Let us see once more how the fBm $W^H\left(u\right)$ appears in this limit
likelihood ratio.
 Denote 
$
\vartheta_u=\vartheta _0+\varphi _Tu,\; g\left(x\right)=  a\, \sgn\left(x \right)\left|x
  \right|^\kappa
$
and write
\begin{align*}
\ln Z_T\left(u\right)&=\int_{0}^{T}\left(S\left(X_t-\vartheta
_u\right)-S\left(X_t-\vartheta _0\right) \right)\:{\rm d}W_t\\
&\qquad \quad -\frac{1}{2}\int_{0}^{T}\left(S\left(X_t-\vartheta
_u\right)-S\left(X_t-\vartheta _0\right) \right)^2\: {\rm d}t\\
& =\int_{0}^{T}\left(g\left(X_t-\vartheta _u\right)-g\left(X_t-\vartheta
_0\right) \right)\:{\rm d}W_t\\
&\qquad \quad -\frac{1}{2}\int_{0}^{T}\left(g\left(X_t-\vartheta
_u\right)-g\left(X_t-\vartheta _0\right)\right)^2\: {\rm d}t+o\left(1\right)
\end{align*}
because
\[
\int_{0}^{T}\left(h\left(X_t-\vartheta _u\right)-h\left(X_t-\vartheta
_0\right)\right)^2\: {\rm d}t\leq C\:u^2 \varphi _T^2\:
T=C\,u^2T^{-\frac{2}{2\kappa +1}+1}\longrightarrow 0 .
\]
Let us denote $\Lambda _T\left(x\right)$ the {\it local time} of the diffusion
process \eqref{4-0} and put $f_T^\circ\left(x\right)=2T^{-1}\Lambda
_T\left(x\right)$. Recall that $f_T^\circ\left(x\right) $ is the {\it local
  time estimator} of the invariant density.  This estimator is consistent
$\bigl($i.e.\ $f_T^\circ\left(x\right)\rightarrow f\left(\vartheta
_0,x\right)\bigr)$ and $\sqrt{T}$-asymptotically normal. Recall that for any
continuous function $H\left(x\right)$ we have
\[
\frac{1}{T}\int_{0}^{T}H\left(X_t\right){\rm d}t=\int_{-\infty }^{\infty
}H\left(x\right)f_T^\circ\left(x\right){\rm d}x\longrightarrow \int_{-\infty
}^{\infty }H\left(x\right)f\left(\vartheta _0,x\right){\rm d}x .
\]
We have
\begin{align}
\label{BT}
\BB_T\left(u\right)&=\frac{a^2}{2}\int_{0}^{T}\left[s_{X_t,u}\left|X_t-\vartheta_u
  \right|^\kappa -s_{X_t,0}\left|X_t-\vartheta_0 \right|^\kappa\right]^2{\rm
  d}t\notag\\
&=\frac{a^2T}{2}\int_{-\infty }^{\infty }\left[s_{x,u}\left|x-\vartheta_u
  \right|^\kappa -s_{x,0}\left|x-\vartheta_0
  \right|^\kappa\right]^2f_T^\circ\left(x\right){\rm d}x\notag\\
&=\frac{a^2T}{2}\int_{-\infty }^{\infty }\left[\sgn\left(y-\varphi
  _Tu\right)\left|y-\varphi _Tu \right|^\kappa -\sgn\left(y\right)\left|y
  \right|^\kappa\right]^2f_T^\circ\left(\vartheta _0+y\right){\rm
  d}y\notag\\
&=\frac{a^2T\varphi _T^{2\kappa +1}}{2}\int_{-\infty }^{\infty
}\left[\sgn\left(v-u\right)\left|v-u \right|^\kappa -\sgn\left(v\right)\left|v
  \right|^\kappa\right]^2f_T^\circ\left(\vartheta _0+\varphi _T v\right){\rm
  d}v\notag\\
&=\frac{a^2 \left|u\right|^{2\kappa +1}}{2}\int_{-\infty }^{\infty
}\left[\sgn\left(s-1\right)\left|s-1 \right|^\kappa -\sgn\left(s\right)\left|s
  \right|^\kappa\right]^2f_T^\circ\left(\vartheta _0+\varphi _T u s\right){\rm
  d}s\notag\\
&\longrightarrow \frac{ \left|u\right|^{2\kappa +1}}{2}a^2\Gamma_*^2
f\left(\vartheta _0,\vartheta _0 \right)=\frac{ \left|u\right|^{2\kappa
    +1}}{2}a^2\Gamma_*^2 G= \frac{ \left|u\right|^{2\kappa +1}}{2}\;\gamma ^2.
\end{align}

Let us denote
\[
\AA_T\left(u\right)=a\int_{0}^{T}\left[\sgn\left(X_t-\vartheta
  _u\right)\left|X_t-\vartheta _u\right|^\kappa -\sgn\left(X_t-\vartheta
  _0\right)\left|X_t-\vartheta _0\right|^\kappa\right]\:{\rm d}W_t.
\]
From the convergence \eqref{BT} and the central limit theorem for stochastic
integrals we obtain
\[
\AA_T\left(u\right)\Longrightarrow \gamma W^H\left(u\right).
\]
More detailed analysis allows verify  the convergence of the finite-dimensional distributions
\[
\Bigl(\AA_T\left(u_1\right),\ldots,\AA_T\left(u_k\right)\Bigr)\Longrightarrow
\Bigl(\gamma W^H\left(u_1\right),\ldots,\gamma W^H\left(u_k\right)\Bigr). 
\]
 Therefore we
have \eqref{4-1}. 

For the MLE $\hat\vartheta _T$ and the BE
$\tilde\vartheta _T$ we have the convergences
\[
T^{\frac{1}{2\kappa +1}} \left(\hat\vartheta _T-\vartheta
_0\right)\Longrightarrow \hat \xi ,\qquad \qquad\qquad \qquad
T^{\frac{1}{2\kappa +1}} \left(\tilde\vartheta _T-\vartheta
_0\right)\Longrightarrow \tilde \xi .
\]
Once more we have the convergence of all polynomial moments and the BE are
asymptotically efficient. 
For the detailed proof see \cite{DK03}, \cite{Kut04}. Note that the case
$\kappa \in \left(-\frac{1}{2},0\right)$ was discussed in \cite{Fu10}.

\subsection{Dynamical system with small noise}

Suppose that the observed process $X^T=\left(X_t,0\leq t\leq T\right)$ is a
solution of the stochastic differential  equation
\[
{\rm d}X_t=\left[a\,\sgn\left(X_t-\vartheta \right)\left|X_t-\vartheta \right|^\kappa
  +h\left(X_t-\vartheta \right)\right]{\rm d}t+\varepsilon {\rm d}W_t, \quad
0\leq t\leq T, 
\]
where the initial value $X_0=x_0 $ is deterministic, $a>0$, $\kappa \in
\left(0,\frac{1}{2}\right)$ and the function $h\left(\cdot \right)$ is known
and has bounded derivative. Moreover we suppose that the function $S\left(x\right)=a\,\sgn\left(x \right)\left|x \right|^\kappa
  +h\left(x \right)>0$ for all $x$. We have to estimate $\vartheta $ and describe the
asymptotic ($\varepsilon \rightarrow 0$) properties of the MLE $\hat\vartheta
_\varepsilon$ and BE $\tilde\vartheta _\varepsilon$.  The likelihood ratio
function is  (see \cite{LS01})
\[
V\left(\vartheta ,X^T\right)=\exp\left\{\int_{0}^{T}
\frac{S\left(X_t-\vartheta\right)}{\varepsilon ^2}{\rm d}X_t-\int_{0}^{T}
\frac{S\left(X_t-\vartheta \right)^2}{2\varepsilon ^2}{\rm d}t\right\},\quad
\vartheta \in\Theta . 
\]
 The set $\Theta $ will be defined below.

We have to  verify the convergence of the normalized likelihood ratio 
\[
Z_\varepsilon \left(u\right)=\frac{V\left(\vartheta_0+\varphi _\varepsilon u
  ,X^T\right) }{V\left(\vartheta_0 ,X^T\right)} \Longrightarrow
Z\left(u\right)=\exp\left\{\gamma
W^H\left(u\right)-\frac{\left|u\right|^{2H}}{2}\gamma ^2\right\}.
\]
Here $\varphi _\varepsilon =\varepsilon ^{\frac{1}{\kappa +\frac{1}{2}}}$ and
$\gamma =a\Gamma_* h\left(0\right)^{-1/2}$.  

The stochastic process $X_t$ converges uniformly on $t\in \left[0,T\right]$ to
$x_t=x_t\left(\vartheta \right) $ --- solution of the ordinary differential equation
\[
\frac{{\rm d}x_t}{{\rm d}t}=a\,\sgn\left(x_t-\vartheta \right)|x_t-\vartheta _0|^\kappa
+h\left(x_t\right),\quad x_0, \quad 0\leq t\leq T.
\]
 Suppose that $\vartheta \in
\left(\alpha ,\beta \right)$, where $ \alpha >x_0$ and $\beta <\inf_{\left\{\alpha
  <\theta\right\} }x_T\left(\vartheta \right)$. 

We can write 
\begin{align*}
\BB_\varepsilon \left(u\right)&=\frac{1}{2\varepsilon ^2 }\int_{0}^{T}
   {\left(S\left(X_t-\vartheta_u \right)-S\left(X_t-\vartheta_0
     \right)\right)^2}{}{\rm d}t \\
&=\frac{a^2}{2\varepsilon ^2}\int_{0}^{T} {\left(s_{X_t,u}\left|X_t-\vartheta
     _u\right|^\kappa - s_{X_t,0}\left|X_t-\vartheta _0\right|^\kappa
     \right)^2}{}{\rm d}t +o\left(1\right)\\
&=\frac{a^2}{2\varepsilon ^2}\int_{0}^{T} {\left(s_{x_t,u}\left|x_t-\vartheta
     _0-\varphi _\varepsilon u\right|^\kappa - s_{x_t,0}\left|x_t-\vartheta
     _0\right|^\kappa \right)^2}{}{\rm d}t +o\left(1\right)\\
&=\frac{a^2}{2\varepsilon ^2}\int_{0}^{T}\frac
   {\left(s_{x_t,u}\left|x_t-\vartheta _0-\varphi _\varepsilon u\right|^\kappa
     - s_{x_t,0}\left|x_t-\vartheta _0\right|^\kappa
     \right)^2}{S\left(x_t-\vartheta _0 \right)}{\rm d}\left(x_t-\vartheta
   _0\right) +o\left(1\right)\\
&=\frac{a^2}{2\varepsilon ^2}\int_{x_0-\vartheta _0}^{x_T-\vartheta _0}\frac
   {\left(\sgn\left(y-\varphi _\varepsilon u\right)\left|y-\varphi
     _\varepsilon u\right|^\kappa - \sgn\left(y\right)\left|y\right|^\kappa
     \right)^2}{ a\, \sgn\left(y\right)\left|y\right|^\kappa
     +h\left(y\right)}{\rm d}y +o\left(1\right)\\
&=\frac{a^2\varphi _\varepsilon ^{2\kappa +1}}{2\varepsilon
     ^2}\int_{-\frac{\vartheta _0-x_0}{\varphi _\varepsilon
   }}^{\frac{x_T-\vartheta _0}{\varphi _\varepsilon }}\frac {\left(
     \sgn\left(v-u\right) \left|v- u\right|^\kappa -
     \sgn\left(v\right)\left|v\right|^\kappa \right)^2}{ a\,
     \sgn\left(v\right)\left|v\right|^\kappa \varphi_\varepsilon ^\kappa
     +h\left(v\varphi_\varepsilon \right)}{\rm d}v +o\left(1\right)\\
& \longrightarrow \frac{a^2}{2h\left(0\right)}\int_{-\infty }^{\infty }
   {\left(\sgn\left(v-u\right)\left|v- u\right|^\kappa -
     \sgn\left(v\right)\left|v\right|^\kappa \right)^2}{\rm
     d}v=\frac{\left|u\right|^{2\kappa +1}}{2}\gamma ^2 .
\end{align*}
Using the same change of variables as above we have
\begin{align*}
\AA_\varepsilon \left(u\right)&=\frac{1}{\varepsilon }\int_{0}^{T}
   {\left(S\left(X_t-\vartheta_u \right)-S\left(X_t-\vartheta_0
     \right)\right)}{\rm d}W_t \\
&=\frac{a}{\varepsilon }\int_{0}^{T}\left[s_{X_t,u}\left|X_t-\vartheta
     _0-\varphi _\varepsilon u \right|^\kappa -s_{X_t,0}\left|X_t-\vartheta
     _0\right|\right]{\rm d}W_t+o\left(1\right)\\
&=\frac{a}{\varepsilon }\int_{0}^{T}\left[s_{x_t,u}\left|x_t-\vartheta
     _0-\varphi _\varepsilon u \right|^\kappa -s_{x_t,0}\left|x_t-\vartheta
     _0\right|\right]{\rm d}W_t+o\left(1\right)\\
&\Longrightarrow \frac{a}{\sqrt{h\left(0\right)}}\int_{-\infty }^{\infty
   }\left[ \sgn\left(v-u\right)\left|v-u\right|^\kappa-
     \sgn\left(v\right)\left|v\right|^\kappa \right]{\rm d}W\left(v\right).
\end{align*}

As above this leads to the following limit distributions for 
the MLE $\hat\vartheta _\varepsilon $ and BE $\tilde\vartheta _\varepsilon $
\[
\varepsilon ^{-\frac{1}{\kappa +\frac{1}{2}}}\left(\hat\vartheta _\varepsilon
-\vartheta _0\right)\Longrightarrow \hat \xi ,\qquad \qquad \varepsilon
^{-\frac{1}{\kappa +\frac{1}{2}}}\left(\tilde\vartheta _\varepsilon -\vartheta
_0\right)\Longrightarrow \tilde \xi ,
\]
convergence of moments of these estimators and the asymptotic efficiency of
the BE. For the full proofs see \cite{K17}.

\section{Discussion}

The proofs presented above can be applied to the
other models of observations. For example, suppose that  we have a nonlinear
stationary time series
\[
X_{j+1}=a\left|X_j-\vartheta \right|^\kappa +h\left(X_j-\vartheta
\right)+\varepsilon _j,\quad j=1,\ldots, n
\]
with i.i.d.\ noise $\left(\varepsilon _j\right)_{j\geq 1}$. Assume the density
$q\left(x\right)$ of the random variable $\varepsilon _j$ and the function $h\left(x\right)>0$ are
sufficiently smooth functions. The likelihood function is
\[
V\left(\vartheta
,X^n\right)=\prod_{j=1}^{n-1}q\left(X_{j+1}-a\left|X_j-\vartheta
\right|^\kappa -h\left(X_j-\vartheta \right)\right),\qquad \vartheta \in
\Theta.
\]
Introduce the notation: $q_{j,u}=q\left(X_{j+1}-a\left|X_j-\vartheta_u \right|^\kappa
-h\left(X_j-\vartheta_u \right)\right)$, $\vartheta _0$ is the true value, $\vartheta_u
=\vartheta_0+\varphi _nu $ and  $g_{j,u}=a\left|X_j-\vartheta_u
\right|^\kappa$. Then  following the same steps as in the i.i.d.\ case above the normalized 
log-likelihood  can be  written as follows  
\begin{align*}
\ln Z_n\left(u\right)&=\sum_{j=1}^{n-1}\ln \frac{q_{j,u}}{q_{j,0}}=
\sum_{j=1}^{n-1}\ln \left[1+\frac{q_{j,u}-q_{j,0}}{q_{j,0}}\right]\\
&=\sum_{j=1}^{n-1} \frac{q_{j,u}-q_{j,0}}{q_{j,0}}
-\frac{1}{2}\sum_{j=1}^{n-1} \left(\frac{q_{j,u}-q_{j,0}}{q_{j,0}}
\right)^2+o\left(1\right)\\
&=\sum_{j=1}^{n-1} \frac{g_{j,u}-g_{j,0}}{q_{j,0}}q'\left(\varepsilon
_j\right) -\frac{1}{2}\sum_{j=1}^{n-1} \left(\frac{g_{j,u}-g_{j,0}}{q_{j,0}}
\right)^2q'\left(\varepsilon _j\right)^2+o\left(1\right)\\
&=a\sum_{j=1}^{n-1} \frac{\left|X_j-\vartheta_0-\varphi _nu
  \right|^\kappa-\left|X_j-\vartheta_0 \right|^\kappa }{q\left(\varepsilon
  _j\right)}q'\left(\varepsilon _j\right)\\
&\qquad - \frac{a^2}{2}\sum_{j=1}^{n-1}
\left(\frac{\left|X_j-\vartheta_0-\varphi _nu
  \right|^\kappa-\left|X_j-\vartheta_0 \right|^\kappa }{q\left(\varepsilon
  _j\right)} \right)^2q'\left(\varepsilon _j\right)^2+o\left(1\right)
\end{align*}
and so on. For nonlinear regression models with cusp-type singularity the
properties of estimators were studied in \cite{PR04}, \cite{D15} and
\cite{DJ15}.

Another interesting problem to discuss is  the estimation of the other parameters of
the model. For example, consider the simplest model 
\[
{\rm d}X_t=a\left|t-b \right|^\kappa {\rm d}t+\varepsilon {\rm d}W_t,\quad
X_0=0,\quad 0\leq t\leq T.
\]
Remind that the parameter $\varepsilon \in \left(0,1\right)$ can be estimated
without error as follows. By It\^o formula for $X_t^2$ we have for any $t\in (0,T]$
\[
X_t^2=2\int_{0}^{t}X_s{\rm d}X_s+\varepsilon ^2 t, \quad {\rm and}\quad
\varepsilon ^2=t^{-1}X_t^2-2t^{-1}\int_{0}^{t}X_s{\rm d}X_s .
\]
The problem of estimation $\vartheta =\left(a,\kappa \right)$ is regular and
the MLE and BE of this parameter are consistent and asymptotically normal with
the regular rate~$\varepsilon$ (see, e.g.\ \cite{IH1}). There is no difficulty
to describe the behavior of the MLE and BE in the case $\vartheta
=\left(a,b\right)$, where $\kappa \in \left(0,\frac{1}{2}\right)$. It can be
shown that, say, the MLE $\hat \vartheta =(\hat a_\varepsilon ,\hat
b_\varepsilon )$ has the following limit distribution
\[
\varepsilon ^{-1}\left(\hat a_\varepsilon -a_0 \right)\Longrightarrow \zeta
\sim {\cal N}\left(0,\sigma ^2\right),\qquad \varepsilon ^{-\frac{1}{\kappa
    +\frac{1}{2}}}\left(\hat b_\varepsilon -b_0 \right)\Longrightarrow
\hat\xi,
\]
where $\zeta $ and $\hat\xi $ are independent random variables.

The problem of estimation
$\vartheta =\left(b,\kappa \right)$ is technically more complicate because the
rate of convergence of the estimator $\hat b_\varepsilon $ depends on the
unknown parameter $\kappa $. It seems that the general results from the
monograph~\cite{IH81} can not be
applied directly here and this problem requires  a special study.

\bigskip

{\bf Acknowledgment.}  This work was done under partial financial support of
the grant RSF number 14-49-00079.  Research of A. Novikov is supported by ARC
Discovery grant DP150102758. The authors are grateful to Dr Lin -Yee Hin for
useful discussions and advices in the simulation study.

\end{document}